\documentstyle{article}

\def\1ox{{ \Omega^1_{\scriptstyle{X}} }}
\def\2ox{{ \Omega^2_{\scriptstyle{X}} }}

\def\ok1{{ \Omega^1_K }}
\def\ok2{{ \Omega^2_K }}

\def\N{{ {\bf N} }}
\def\Q{{ {\bf Q} }}

\def\C{{ {\bf C} }}
\def\R{{ {\bf R} }}
\def\8{{ {\infty } }}

\def\deg{{ \mbox{deg} }}
\def\^{{ ^{\wedge} }}

\def\Z{{ {\bf Z } }}

\def\sup{{ \mbox{sup} }}

\title{Diophantine Equations in Two Variables}
\author{Minhyong Kim}
\begin{document}
\maketitle

We start with an example:

$$x^3+y^3=1729$$

Let's try to find all integral solutions to this
equation. One factors the
polynomial on the left
$$x^3+y^3=(x+y)(x^2-xy+y^2).$$
So any integral solution has to satisfy
$$x^2-xy+y^2\leq 1729$$
hence,
$$(x-y/2)^2+3y^2/4\leq 1729$$
which gives us the bound
$$|y|\leq 2\sqrt{1729/3}\approx 48$$
and similarly for $|x|$. Thus, we can carry out an
exhaustive search for the solutions and find the solutions
$(9,10)$,$(1,12)$ and their reflections. 
This is a deceptively simple example of a Diophantine
problem.
Nevertheless, from the point of view of
complexity, this problem already is non-trivial. That is, if one considers
the family of equations
$$x^3+y^3=m$$
then the same method will give bounds
$$\sup \{|x|,|y|\} \leq 2\sqrt{m/3}$$
using which we can apply the search algorithm.
However, the search time is clearly exponential
in the length of $m$.
So we are saying actually that the problem is very
simple {\em theoretically}. Still, one point
I hope to make clear in this lecture is that
when it comes to Diophantine equations,
one is extremely fortunate even to have
an exponential algorithm, or any algorithm
at all.

 In any case, one can see how deceptive this simplicity
is by considering
the same problem where we ask for all $(x,y)$ in
$\Z[\sqrt{-1}]$ or $\Z[\sqrt{2}]$. Is the solution
set in fact finite in these cases? (Yes.)
The problem is also considerably more difficult if
we consider $\Q$-solutions. For example, one checks
that  $$P=(\frac{-3457}{1727}, \frac{20760}{1727})$$
and $$Q=(\frac{-5150812031}{1075576681},\frac{5177701439}{107557668})$$
   are solutions. How does one find  solutions
involving such large numbers? The answer is that one
applies {\em geometric constructions}. The point $P$
is the intersection between our curve and
the tangent line at $(1,12)$ while $Q$ is the
curve's intersection with the secant line connecting
$P$ and $(1,12)$. It is an easy exercise to
show that such constructions necessarily yield
$\Q$-points.
It is also not hard to show that one obtains in fact infinitely
many $\Q$-points in this manner. (Repeatedly applying
tangents will yield rational solutions with
smaller and smaller denominators.)

The point of this discussion is that 
determining the structure of the solutions
coming from different rings can be
quite simple or complicated depending on the ring.

On the other hand one normally regards the problem of
finding solutions in $\R$ or $\C$ as easy, but of
course this is because one isn't asking for the
same kind of description as in the $\Z$ or $\Q$
case, 
nor would such a description even make
much sense. 

The equations of the title are of the
form
$$X:\ \ \ f(x,y)=0$$
where $f\in R[x,y]$ is a polynomial of two
variables with coefficients in some (small) ring $R$
and the associated problem is that of finding
the structure of the solutions $(x,y)$
where $(x,y)$ are subject to the condition
that they also come from some pre-specified
ring $k$ (possibly large). We will denote such a set
of solutions by $X(k)$ and often refer to
it as the set of $k$-points of $X$.
Any such problem is a Diophantine
problem in two variables. 

A brief remark on $R$ and $k$: for the most part,
you should think of them as subrings of $\C$.
How would they be given? Usually, in terms of generators,
so that we might have
$R=\Z[a_i]$ where the $a_i$ are the roots
of $x^3-x+1$ and $k=\Z[a_i,\pi]$.

 We will be
assuming in this lecture that the polynomial $f$ is geometrically irreducible,
that is, can't be non-trivially factored over $\C$.
Also, just to keep the discussion stream-lined, we
will also assume that $X$ is smooth, that is,
the gradient $\nabla f$ is non-vanishing on
all the complex points of $X$.

Eventually, one would like to understand
the relation between the various sets
$X(k)$ as we vary the ring $k$.
For example, we might start with $k=\Q$,
enlarge it to some field of algebraic numbers,
$k_1$, so that $X(k)\subset X(k_1)$,
and keep going 
$$X(k)\subset X(k_1) \subset X(k_2) \subset \cdots \subset X(\C).$$
One should be aware that there are some incredible
leaps along the way which we are at present
unable to understand in any  satisfactory fashion.
At the end of this increasing filtration, one
gives up on precise arithmetic understanding as the
set acquires a certain homogeneity and is usually
studied as a geometric object. 
At the beginning, set-theoretic considerations
dominate.
From the  point of view of this filtration,
therefore,  the goal of Diophantine
geometry is to give 
descriptions and even constructions of the various
middle stages, or in
other words, to understand the {\em arithmetic formation}
of geometric objects. 
This a grand problem, whose correct formulation
 even is unclear. So we spend most of our time at
the lowest level, trying to understand sets like
$X(\Z)$ or $X(\Q)$, or if we are
really willing to work hard, $X(\Q(\sqrt{-1}))$,
or at the top level, where we study the
various properties of complex algebraic curves.

These simplest
cases are usually enough to occupy
the careers of large collections of number-theorists
and geometers.
Nevertheless, the mainstream of
work in   Diophantine geometry of this
century has clearly demonstrated  at least that
a link exists between the bottom and top
levels, and that
 the coarse
structure of the solution set in large rings
 is related to the coarse structure of the solution
set in small rings. We will see examples of this
phenomenon presently.

Why two variables? The answer is that, unfortunately, this
 is essentially the only situation where we have
something close to a theory. It has to
do with one of the basic invariants of a Diophantine
system, namely, the dimension. At present, a coherent
discussion needs to limit itself to the
case of dimension 1, which, more or less, is
the study of two-variable equations.

As mentioned the basic questions at the bottom level are
essentially set-theoretic:

Is $X(k)$ non-empty?

Is it finite or infinite?

Can we determine the set $X(k)$?

We go on to discuss some of the answers which flesh
out the aforementioned theory of Diophantine
equations in two variables. For simplicity,
we will assume that $f(x,y)\in \Z[x,y]$ and
for the most part discuss only $X(\Z)$ and $X(\Q)$.
First we need to recall one other invariant of
the equation (*). This is its genus, which we
will denote $g=g(f)$, and it
can be computed as follows: One considers
the set $X(\C)$ of complex solutions to the equation.
This is a topological space of real dimension two.
One can show (using projective closure and
resolution of singularities) that this is in
fact a smooth compact surface with a few
points removed. The genus of $f$ is by definition  the
genus of this compact surface. (Incidentally
algebraic geometers are more fond of calling this
object a curve in view of its complex dimension.)
It is easy to compute that generically  $g=0$ for degree 1 and 2,
$g=1$ for degree 3, and $g\geq 2$ for $d\geq 4$.
The fact that this simple formula is true only
generically is actually rather important,
and it was a significant realization in Diophantine
geometry that the  genus is in fact a more
basic invariant than the degree.

Note that this invariant of the equation is read off at the
top level. Nevertheless, it influences the
set-theoretic  structure at the
bottom level  via the following trichotomy:
$$\begin{array}{cc}
g=0 & X(\Q)=\phi \ \ \  \mbox{or} \ \ \ \infty \\
g=1 & X(\Q)=\phi ,\ \ \  \mbox{finite}\ \ \  \mbox{or} \ \ \  \infty \\
g\geq 2 & X(\Q)= \phi \ \ \  \mbox{or} \ \ \  \mbox{finite}\\
\end{array}$$

This list should be supplemented by the fact that the genus
1 case has much less points than the genus zero case in
a very precise sense, even when both are infinite.
The hardest fact is in the third line, which was a
conjecture of Mordell, resolved by Faltings in 1983.
In fact, we could go on to turn the last line
into a characterization as follows:
The genus of $f$ is $\geq 2$ iff $X(k)$ is
finite for every finitely generated field $k$.

There is also a corresponding statement for
integral points. For this, we need 
a symbol $s$ for the number of points removed
from the compact surface to get $X(\C)$. (This is
usually equal to the number of linear factors
in the higheest degree homogenous part of $f(x,y)$.)

The theorem of Siegel says then that
$2g-2+s>0$ iff $X(R)$ is finite for
every  finitely generated ring $R$.
It is interesting to formulate
the condition $2g-2+s>0$ geometrically:
$X(R)$ is finite for every finitely generated
ring $R$ iff $X(\C)$ admits a hyperbolic metric.

In any case, one catches a glimpse of
the deep relations that can exist between
$X(k)$ for different rings $k$.

Let us try to get a better feeling for the
trichotomy through some examples,
starting with the case of genus zero.
A reasonably general class of
genus one equations is the family
$$ax^2+by^2=c$$
where $a,b,c$ are positive integers.
There is a standard procedure for geometrically generating
all rational solutions out of one which is best
illustrated by the example of the circle
$x^2+y^2=1$
One shows easily that all rational points
are obtained by intersecting the
circle with the lines going through the point $(-1,0)$
with rational slope $t$. In fact, the point corresponding
to slope $t$ has the explicit formula
$$(\frac{1-t^2}{1+t^2},\frac{2t}{1+t^2})$$
and I leave it to the audience to convince themselves
that the same technique will hold for the general
equation in our family. The only problem in the general
case is that the single point which we need to
get going might not be evident. In fact, it may not exist
as we see from the trivial example
$x^2+y^2=-1$. But how about $x^2+y^2=3$?
An algorithm for finding a solution if it exists
is provided by a nice theorem of Holzer:
Suppose the equation $ax^2+by^2=c$ has a rational
solution. Then there is a solution $(p/r,q/r)$ (in reduced
form) with $\sup(|p|,|q|,|r|)\leq \sqrt{abc}$. This
immediately provides verification that
$x^2+y^2=3$ has no solutions.
Notice that one can view this theorem as providing
a decision procedure to the question of whether
or not the solutions set is empty. However, this is
again an exponential-time algorithm.
If one uses the Hasse principle, it is possible to
get a subexponential algorithm, that is, one that
requires determining the prime factors of $abc$, but
still not polynomial.
It is interesting to note that the method for
finding rational solutions a priori says nothing
about the existence of integral solutions.
The latter are found by viewing the equation
as a norm equation and using techniques of
algebraic number theory. The process 
is somewhat more involved than the case of rational solutions
and we will bypass it in this lecture.
Incidentally, if one changes the equation slightly,
a theorem of Adleman and Manders states that
determining the existence of natural number solutions to
$ax^2+by=c$ is NP-complete.

Genus one equations that have at least one solution
can more or less be transformed into
the form
$$E:y^2=x^3+ax+b$$
(exercise: try this for the
first example of this lecture)
and these are among the most studied equations
in number theory.  As mentioned, the set
of integral solutions is finite, while the
rational solutions can be infinite.
Nevertheless, there is a remarkable
fact due to Mordell that says
that $E(\Q)$ can be generated
by a finite set of solutions using the geometric
process of constructing secants and tangents.
In this respect, the situation is analogous
to the genus zero case. That is,
even though we are trying to describe an infinite
set, it is in essence `contained' inside
its finite set of generators.
 However, finding a generating
set turns out to be much more difficult than
for genus 0. In fact, we have at present no
algorithm for carrying this out. The
existence of a good algorithm is an important
part of the conjecture of Birch and Swinnerton-Dyer.
That is, the finiteness of an invariant called
the Shafarevich group will give us in principle
an algorithm for finding a generating set for
$E(\Q)$.  However, even in the best case the complexity
appears to be quite high, certainly at least exponential,
and I have been too lazy
to examine it prior to this lecture.
There is a well-known
 algorithm called the method of two-descent
which works well when the 2-part of the Tate Shafarevich group
is trivial and it has been implemented by John Cremona.
We note that it is also possible to give a conditional
algorithm based upon the $L$-function of
the elliptic curve again assuming the
B-Sw-D conjecture.
In fact, this forms the
basis of quite a bit of activity surrounding
integral solutions. The fact that these could all be
found in principle arose from the celebrated
work of Alan Baker which implied that all
integral solutions $(x,y)$ must satisfy
the bound:
$$\sup \{|x|,|y|\} \leq \exp ((10^6 \sup \{|a|,|b|\})^{10^6})$$
Of course, this bound would be completely
useless in practice, giving us only
doubly exponential complexity.
However, ideas of
Lang and Zagier on the elliptic logarithm,
transcendence results of David, and
the $L^3$ algorithm for lattices have been
combined to devise an algorithm
that can often be used to find solutions in exponential  time,
and even polynomial time in fortunate circumstances.
Using these methods, various established cases
of the B-Sw-D conjecture, and occasional triviality
of the Shafarevich group, Gebel, Peth\"{o}, and Zimmer
have succeeded in finding all integral points
on
$$y^2=x^3+k$$
for all $|k|\leq 10000$ and most $|k|\leq 100000$.
From outside the field, I suspect that it is
quite surprising to hear that these computations
are as recent at 1998 and that this is not for lack
of trying. Also astounding is the triviality of the
achievement from a complexity-theoretic standpoint,
involving merely 5 or 6 digit inputs.

When we move to higher genus curves, which
we get for the general polynomial $f(x,y)$
of degree at least 4,
even though we
know the rational solution set to be
always finite, it is ironical that
one typically needs to work very hard
and have a good deal of luck to deal with
any given equation (recall the case
of the Fermat equation). Feasible methods  for
computing integral points exist for very special
families such as Thue curves (due to De Weger
and Tzanakis), and for rational
points the method of Coleman-Chabauty has had
some success in provably finding all rational points
on certain hyperelliptic curves.

Rather than survey those results, we will spend the rest of our time
leading up to the so-called {\em effective Mordell
conjecture}.

By way of motivation, let us digress a moment
to discuss the case of Diophantine
equations over function fields.
That is to say, in the language we
introduced
at the beginning, we are interested in the
case where $R=k=\C(t)$. 

Example:
$$(\frac{1-t^2}{1+t^2})^2+(\frac{2t}{1+t^2})^2=1$$

More generally, we are studying equations like
$$\sum_{i,j}a_{ij}(t)x^iy^j=0$$

It is probably well-known that
 techniques of complex geometry  (the
study of algebraic fibrations) render these
problems much more tractable than their arithmetic
versions. For example, just complex
analysis will suffice to show that
the Fermat equation
$x^n+y^n=1$
for $n\geq 3$ has no non-constant
solutions in rational functions.
In general, if we have an integral coefficient
equation so that we can consider both
rational number and non-constant rational function solutions,
it is easy to see  from topological considerations
that the latter are much hard to come by.

It is perhaps less known that the finiteness
theorems for higher genus equations
over function fields often come with
effective a priori estimates for
the complexity of the solutions.
For example,
Let $f(t,x,y)\in \C[t,x,y]$ be a generic polynomial of
$x,y$ degree at least four, which we are therefore
viewing as a polynomial in $x,y$ with coefficients in
$\C[t]$. (Here we are merely using the term generic to
keep the statements simple, and one can actually
deal with the general case if one is willing to divide
into many subcases, or use algebraic-geometric
terminology.) Then
$$f(t,x,y)=0$$
has only finitely many solutions
in rational functions and any solution
 $$(p(t)/r(t), q(t)/r(t))$$where $(p,q,r)=1$
satisfies the bound
$$\sup \deg \{p,q,r \} \leq (d^2-3d-1)(2s+1) $$
where $s$ is the number of zeroes of the polynomial in
$t$ obtained by computing the $(x,y)$ discriminant of $f$.
Actually, once one has this bound, it is possible in principle
to compute all solutions simply by substituting rational
functions with undetermined coefficients and using a
Groebner basis algorithm. 

So here is the effective Mordell conjecture for
Diophantine equations in two variables:
Given a polynomial $f(x,y)\in \Z[x,y]$ of genus at least
two, there is a computable constant $C(f)$
such that all rational number solutions $(p/r,q/r)$ with $(p,q,r)=1$
to the equation
$$f(x,y)=0$$
satisfy the bound $\sup\{|p|,|q|,|r|\}\leq C(f)$.
Another way to view this is to simply consider the sup over
all points of the sup of numerators and denominators. 
That is, define
$$A(f)=\sup_{(p/r,q/r)\in X(\Q)} \sup\{|p|,|q|,|r|\}$$
This clearly exists since we know the
solutions set to be finite and it can be viewed just
as a function of $f$. The conjecture then is that $A(f)$
is actually computable. This is a very bold conjecture in
that one knows the Diophantine decision problem to be
undecidable even for equations in 9 variables.
It should be emphasized that any bound at all
for $A(f)$, even as big as
one of Harvey Friedman's numbers would be a remarkable theorem. 
On the other hand, optimistic people also conjecture
bounds that are not very different from the
 ones obtained over function fields.

Many interesting ideas exist for attacking this
conjecture and there is also a body of related
conjectures all more or less equivalent whose inner
coherence  appears to render them together plausible.

The ideas are mostly inspired by algebraic geometry,
namely, intersection theory and deformation theory.
But it is fair to say at present that there are no
concrete results. One of my reasons for attending
this workshop was  the hope that inspiration from complexity
theory might lead to
a fresh line of attack.

At a workshop of this nature, it might not be entirely
 out of place
to close this lecture with some philosophical
observations. There is a well-known reason for the
difficulty of Diophantine problems: A theorem
of Matiyasevich states that almost any class of
mathematical problems can be reduced to non-existence
of positive integral solutions to a Diophantine equation.
For example, there is a polynomial $P(x_1,\ldots, x_n)\in \Z [x_1,\ldots, x_n]$
such that the Riemann hypothesis is true iff
$$P(x_1,\ldots, x_n)=0$$
has no solution. One can similarly find an equation for, say,
the  four-color problem. We are even quite close to having such
an equation for the Poincare conjecture. 
It should then come as no surprise that
we can relate our main conjecture
$P\neq NP$ to a Diophantine equation, and in a more interesting
way than just cooking up a Diophantine decision problem
that is NP complete.
So if you really knew your Diophantine equations, millions
would be available to you immediately!

Let me briefly outline the  connection.
Suppose $S\subset \N^k$ is a recursively enumerable set.
 Matiyasevich's representation theorem says
that there is a polynomial $f(t_1,\ldots, t_k, x_1,\ldots, x_n)$
with integer coefficients such that
$$f(t_1,\ldots, t_k, x_1,\ldots, x_n)=0$$
has a natural number solution in the $x_i$'s iff $(t_1,\ldots, t_k)\in S$.

Let's recall from Avi's lecture the following equivalent formulation
of $P\neq NP$:
Consider ways of writing 
$$\mbox{perm}_n =\mbox{det}_m\circ A$$
where $A$ is an affine linear map and we consider this as an
identity of polynomials over ${\bf F}_2$. 
Such a representation exists for $m=2^n$. Finally,
recall  that $P\neq NP$ is equivalent to
the assertion that the minimal $m$ for which such a representation
exists is not polynomial in $n$. Now, for a given $n$, we
can clearly compute the minimal $m$ via enumeration.
Call this minimum $l(n)$. Then $P\neq NP$ says
that there are no natural numbers $a$ and $k$
such that $$l(n) \leq an^k$$ for all $n$. But
the set $S$ of $(a,k)$ for which some $n$ satisfies
$l(n)> an^k$ is recursively enumerable. (Order
$(a,k,n)$ according to some non-decreasing function
of $a+k+n$ and go through them one by one. Toss $(a,k)$
into our set $S$ whenever $l(n)>an^k$.)
Now, let $f(t_1,t_2, x_1,\ldots, x_n)$ be
the polynomial from Matiyasevich's theorem.
Then $P\neq NP $ is equivalent to
the assertion that 
$f(t_1,t_2, x_1,\ldots, x_n)=0$
has a solution for all $(t_1,t_2)$. To appreciate this
formulation, note that if you could show for any
fixed $(a,k)$
that the Diophantine equation
$f(a,k, x_1,\ldots, x_n)=0$ 
has {\em no} solution, then you would have proved in fact that
$P=NP$!

Although there are undoubtedly some serious aspects
to this line of thought, as I'm sure I've exhausted
the hour and your patience, this seems a good place to
stop.

\end{document}